\documentclass[11pt]{article}
\usepackage[margin=1in]{geometry}
\usepackage[T1]{fontenc}
\usepackage{lmodern,microtype,amsmath,amssymb,amsthm,mathtools}
\usepackage{booktabs,enumitem,tikz,placeins,graphicx}
\usetikzlibrary{arrows.meta}
\usepackage[colorlinks=true,linkcolor=blue!45!black,citecolor=blue!45!black,urlcolor=blue!45!black]{hyperref}
\hypersetup{pdftitle={Improving the Round Taylor Method},pdfauthor={Oscar M. Perdomo}}
\setlist[enumerate]{itemsep=3pt}
\newtheorem{theorem}{Theorem}[section]
\newtheorem{lemma}[theorem]{Lemma}
\theoremstyle{remark}
\newtheorem{remark}[theorem]{Remark}
\newcommand{\R}{\mathbb R}
\newcommand{\norm}[1]{\left\lVert #1\right\rVert_\infty}
\newcommand{\dd}{\,d}
\definecolor{inkblue}{RGB}{28,86,126}
\definecolor{leaf}{RGB}{36,116,87}
\title{Improving the Round Taylor Method and an automatic procedure to solve ODEs with polynomial vector fields }
\author{Oscar M. Perdomo}
\begin{document}
\maketitle

\begin{abstract}
Validated integration already provides methods for certifying
existence and bounding the error of numerical approximations to
solutions of ordinary differential equations. 
 We present a variation of the Round Taylor method, introduced
by the author in 2019, with a compact formulation that is
straightforward to implement. The method certifies existence of the solution throughout the
computed time interval and provides a rigorous global error
bound that accounts explicitly for the errors introduced by
rounding the computed approximations at each step. For polynomial vector fields with rational coefficients, we give
a self-contained presentation of standard interval-range bounds specialized to our setting,
and an automatic validation procedure using exact rational
arithmetic. This implementation does not require a specialized
interval-arithmetic library. Given a target time, the procedure
selects step sizes and computes error bounds, either certifying
an approximation at that time or reporting the step at which
certification stops. We include examples in one and two variables
that can be followed by hand, a scalar computer example with certified
error below $10^{-20}$, and a four-body problem. For the four-body
example, we certify small residuals that give us evidence of the existence
of a periodic motion.
\end{abstract}

\section{Introduction}
Let us start by explaining why rounding can become necessary in
a practical computation with rational numbers. Rational arithmetic
allows each operation to be performed exactly, but the numerators
and denominators can grow rapidly as the computation proceeds.
In a long integration,
storing and manipulating these integers can make the unrounded
computation impractical. This difficulty motivated the rounding
procedure in the author's earlier work on periodic solutions of
the three-body problem~\cite{Perdomo2015}.

Rounding keeps the stored numbers manageable, but it introduces
an error at every step. That error is carried into later steps
and may be amplified. An estimate of the Taylor remainder alone
therefore does not bound the error in the values actually stored.
The numerical method must also account for the initial approximation error
and for all subsequent rounding errors. The Round Taylor
method~\cite{Perdomo2019} does this using rational arithmetic
and prescribed rounding grids. Under explicit hypotheses, its
error estimate also guarantees existence of the exact solution
throughout the computed time interval.

These goals and several of the validation principles used here
are established in the literature on validated integration.
Nedialkov, Jackson, and Corliss~\cite{NJC1999} survey
Taylor-series methods for validated initial-value problems.
Jackson and Nedialkov~\cite[Sections~2--3]{JN2002} describe
interval range bounds, a larger region for validating existence
and the Taylor remainder, and Jacobian-based propagation on the
smaller incoming enclosure. Their work also accounts for
computer rounding. Kapela, Mrozek, Wilczak, and
Zgliczy\'nski~\cite{CAPD2021} describe the CAPD::DynSys toolbox
for rigorous integration and computer-assisted proofs in dynamics.
Our presentation is informed by this existing framework and
is intended as a compact alternative formulation that requires no specialized interval-arithmetic library.

Our new version of the Round Taylor theorem replaces the
separate derivative bounds in the 2019 formulation by one bound
for the derivative of the complete Taylor predictor. It also
allows the certifying regions and bounds to change from step to
step. The incoming-error cube is used for the predictor derivative,
and a larger neighborhood is used for the speed and Taylor
remainder. This formulation makes the sufficient conditions
explicit while retaining prescribed rounding grids.

For polynomial vector fields, the bounds follow standard
interval-range principles, such as those described
in~\cite[Section~2.1]{JN2002}. Our polynomial section provides a
self-contained specialization and an exact rational
implementation of these bounds. The resulting automatic
procedure requires no specialized interval-arithmetic library.
It computes the steps and required bounds and returns either a
certified approximation at the target time or the last certified
point and a reason for stopping. The hand examples follow this
procedure in the order in which the computations are performed.
The computer examples use moderate Taylor orders and state the
rounding grid and the resulting error explicitly. The four-body
example also explains how to remove radicals from the vector field
and enclose the initial reciprocal distances~\cite{Perdomo2026NBody}.
We then compute to the symmetry midpoint used in that revised
manuscript and bound the residuals of its four symmetry conditions.
We explain why exact satisfaction of these conditions would give
a periodic motion.

\section{Preliminaries and the definition of the method}
For vectors and matrices, respectively, use
\[
 \norm v=\max_{1\le i\le n}|v_i|,\qquad
 \norm M=\max_{1\le i\le n}\sum_{k=1}^n|m_{ik}|.
\]
Then $\norm{Mv}\le\norm M\norm v$. If $g$ is continuously
differentiable and $\norm{Dg}\le A$ on the segment from $a$ to
$b$, integration along that segment gives
\begin{equation}\label{eq:meanvalue}
 \norm{g(b)-g(a)}\le A\norm{b-a}.
\end{equation}

Consider $Y'=f(Y)$ with initial condition $Y(0)=y_0$, where
$f\in C^{m+1}(\Omega,\R^n)$, $\Omega$ is open, and $m\ge1$.
The Taylor functions are
\[
 F_0=f,\qquad F_{\ell+1}=DF_\ell f
 \quad(\ell=0,\ldots,m-1).
\]
Along a local exact solution, $Y^{(\ell+1)}=F_\ell(Y)$ for
$\ell=0,\ldots,m$.

For a positive rational grid spacing $H$, define
\[
 {\cal R}_H(v)=H\left\lfloor\frac vH\right\rfloor,
\]
with the floor applied coordinatewise. For every $v\in\R^n$,
the vector ${\cal R}_H(v)$ is rational because $H$ is rational
and each floor is an integer. Since
$0\le a-H\lfloor a/H\rfloor<H$ for every real $a$,
we have $\norm{v-{\cal R}_H(v)}<H$.

For $N$ steps, choose positive rational numbers
$h_0,\ldots,h_{N-1}$ and $H_0,\ldots,H_N$.
Let $y_0\in\Omega$ and set $z_0={\cal R}_{H_0}(y_0)$ and
$E_0=H_0$.
For $j=0,\ldots,N-1$, define
\begin{align}
 t_0&=0,\qquad t_{j+1}=t_j+h_j,\notag\\
 P_{j+1}(y)&=y+\sum_{\ell=1}^m
                \frac{h_j^\ell}{\ell!}F_{\ell-1}(y),
                \label{eq:predictor}\\
 z_{j+1}&={\cal R}_{H_{j+1}}(P_{j+1}(z_j)).
                \label{eq:roundstep}
\end{align}
Thus there are $N$ predictors $P_1,\ldots,P_N$ and $N+1$
stored points $z_0,\ldots,z_N$.

Let $E_1,\ldots,E_N$ be proposed subsequent error bounds. At each step
$j=0,\ldots,N-1$, we use a speed bound $B_j$ and a Taylor
remainder bound $\tau_j$. A speed bound explains the larger
radius $\varepsilon_j$: while the solution remains in a region
where $\norm f\le B_j$, its distance from $z_j$ is at most
$E_j+B_js$ after time $s$. We therefore require
$E_j+B_jh_j<\varepsilon_j$.

The derivative bound $A_j$ compares two predictor values,
while $\eta_j$ bounds the numerical defect. The corresponding
error inequality is
\[
 E_{j+1}\ge A_jE_j+\tau_j+\eta_j.
\]
The derivative is bounded on the incoming-error cube $S_j$;
the speed and remainder are bounded on the larger cube $C_j$.
The following theorem makes these requirements precise.

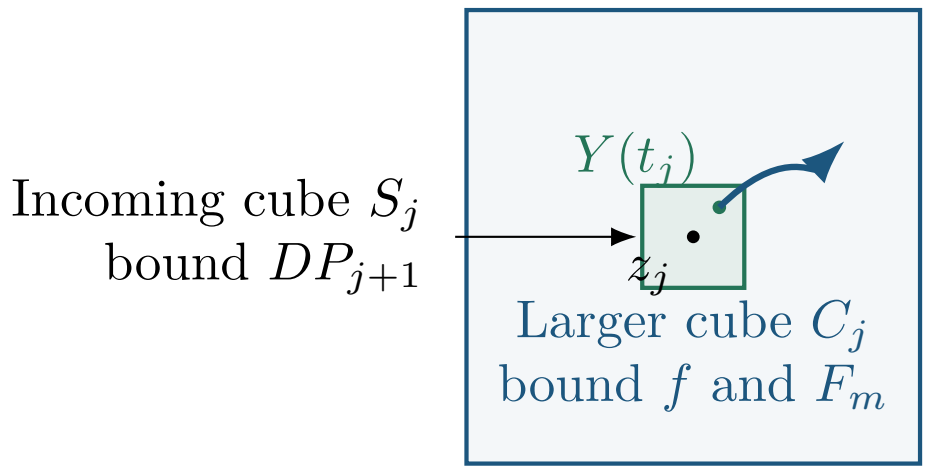
\begin{figure}[ht]\centering
\begin{tikzpicture}[scale=.75,>=Latex,font=\small]
\fill[inkblue!5] (-2,-2) rectangle (2,2);
\draw[inkblue,thick] (-2,-2) rectangle (2,2);
\fill[leaf!12] (-.45,-.45) rectangle (.45,.45);
\draw[leaf,thick] (-.45,-.45) rectangle (.45,.45);
\fill (0,0) circle (1.6pt) node[below left]{$z_j$};
\fill[leaf] (.23,.26) circle (1.7pt) node[above left]{$Y(t_j)$};
\draw[inkblue,very thick,->] (.23,.26) .. controls (.75,.78) and (1.05,.57) .. (1.34,.85);
\node[anchor=east,align=right] at (-2.2,0)
 {Incoming cube $S_j$\\bound $DP_{j+1}$};
\draw[->] (-2.1,0)--(-.5,0);
\node[inkblue,align=center] at (0,-1.05)
 {Larger cube $C_j$\\bound $f$ and $F_m$};
\end{tikzpicture}
\caption{The incoming state lies in $S_j$; the inequality
$E_j+B_jh_j<\varepsilon_j$ keeps the exact motion inside $C_j$.}
\label{fig:neighborhoods}
\end{figure}
\FloatBarrier

\section{Main result}\label{sec:theorem}

\begin{theorem}[improved Round Taylor]\label{thm:main}
Let $n,m,N$ be positive integers, let $\Omega\subset\R^n$ be open,
and let $f\in C^{m+1}(\Omega,\R^n)$. Consider the initial-value problem
\[
 Y'(t)=f(Y(t)),\qquad Y(0)=y_0\in\Omega.
\]

Let $h_0,\ldots,h_{N-1}$ and $H_0,\ldots,H_N$ be positive
rational numbers. Set $z_0={\cal R}_{H_0}(y_0)$ and $E_0=H_0$,
and suppose that the following sequences are given:
\begin{itemize}
\item nonnegative error bounds $E_1,\ldots,E_N$
      and positive radii $\varepsilon_0,\ldots,\varepsilon_{N-1}$;
\item nonnegative speed bounds $B_0,\ldots,B_{N-1}$
      and derivative bounds $A_0,\ldots,A_{N-1}$;
\item nonnegative remainder bounds $\tau_0,\ldots,\tau_{N-1}$
      and numerical-defect bounds $\eta_0,\ldots,\eta_{N-1}$.
\end{itemize}

Define the Taylor functions recursively by
\[
 F_0=f,\qquad F_{\ell+1}(y)=DF_\ell(y)f(y)
 \quad(\ell=0,\ldots,m-1).
\]

Set $t_0=0$ and, for $j=0,\ldots,N-1$, define
\[
 t_{j+1}=t_j+h_j,\qquad
 P_{j+1}(y)=y+\sum_{\ell=1}^{m}
                \frac{h_j^\ell}{\ell!}F_{\ell-1}(y).
\]

Starting from $z_0$, define the remaining stored points by
\[
 z_{j+1}={\cal R}_{H_{j+1}}(P_{j+1}(z_j))
 \quad(j=0,\ldots,N-1).
\]

For each $j=0,\ldots,N-1$, define the two neighborhoods
\[
 S_j=\{y\in\R^n:\norm{y-z_j}\le E_j\},\qquad
 C_j=\{y\in\R^n:\norm{y-z_j}<\varepsilon_j\}.
\]
Assume that $\overline{C_j}\subset\Omega$ for every
$j=0,\ldots,N-1$, and that the following bounds hold for
$j=0,\ldots,N-1$:
\begin{align}
 \norm{f(y)}&\le B_j &&(y\in C_j),\label{eq:speed}\\
 \norm{\frac{h_j^{m+1}}{(m+1)!}F_m(y)}
     &\le\tau_j &&(y\in C_j),\label{eq:remainderbound}\\
 \norm{DP_{j+1}(y)}&\le A_j &&(y\in S_j),\label{eq:mapbound}\\
 \norm{z_{j+1}-P_{j+1}(z_j)}&\le\eta_j.\label{eq:defect}
\end{align}

For every $j=0,\ldots,N-1$, assume also
\begin{align}
 E_j+B_jh_j&<\varepsilon_j,\label{eq:escape}\\
 E_{j+1}&\ge A_jE_j+\tau_j+\eta_j.\label{eq:recurrence}
\end{align}
Then the initial-value problem has a unique solution on $[0,t_N]$,
and
\begin{equation}\label{eq:nodeerror}
 \boxed{\norm{Y(t_j)-z_j}\le E_j\quad(j=0,\ldots,N).}
\end{equation}
During every step it satisfies
\begin{equation}\label{eq:stepcontainment}
 \norm{Y(t_j+s)-z_j}\le E_j+B_js<\varepsilon_j
 \quad(j=0,\ldots,N-1,\quad 0\le s\le h_j).
\end{equation}
\end{theorem}

\begin{proof}
By local existence and uniqueness, the initial-value problem has a
unique maximal solution $Y:[0,T_{\max})\to\Omega$, where
$T_{\max}>0$. We prove by induction that $t_j<T_{\max}$ and
$\norm{Y(t_j)-z_j}\le E_j$ for $j=0,\ldots,N$.
For $j=0$, the rounding bound gives
\[
 \norm{Y(0)-z_0}=\norm{y_0-{\cal R}_{H_0}(y_0)}<H_0=E_0.
\]

Suppose it holds at some $j<N$, and put $X_j=Y(t_j)$.
By~\eqref{eq:escape}, $E_j<\varepsilon_j$, so that
$X_j\in S_j\subset C_j$. As long as the solution remains
in $C_j$, integration of the differential equation gives
\begin{align}
 \norm{Y(t_j+s)-z_j}
 &\le\norm{X_j-z_j}
       +\int_0^s\norm{f(Y(t_j+v))}\dd v\notag\\
 &\le E_j+B_js.\label{eq:displacementproof}
\end{align}
If the solution had a first exit from $C_j$ at a time
$t_j+s_*<T_{\max}$ with $0<s_*\le h_j$, continuity would give
$\norm{Y(t_j+s_*)-z_j}=\varepsilon_j$. Taking a limit
in~\eqref{eq:displacementproof} would instead give
\[
 \norm{Y(t_j+s_*)-z_j}
 \le E_j+B_js_*\le E_j+B_jh_j<\varepsilon_j,
\]
a contradiction. Thus, for every $0\le s\le h_j$ with
$t_j+s<T_{\max}$, the solution lies in the compact set
\[
 K_j=\{y\in\R^n:\norm{y-z_j}\le E_j+B_jh_j\}
 \subset C_j\subset\Omega.
\]

If $T_{\max}\le t_{j+1}$, the speed bound gives
\[
 \norm{Y(t)-Y(s)}\le B_j|t-s|
 \qquad(t_j\le s,t<T_{\max}).
\]
Consequently $Y(t)$ has a limit $X_*\in K_j$ as
$t\uparrow T_{\max}$. Since $X_*\in\Omega$, local existence
at $X_*$ extends $Y$ beyond $T_{\max}$, contradicting its
maximality. Therefore $t_{j+1}<T_{\max}$, and
\eqref{eq:displacementproof} holds for every $0\le s\le h_j$.
In particular, the solution remains in $C_j$ throughout this step.

The chain rule gives $Y^{(m+1)}(t)=F_m(Y(t))$.
Taylor's formula with integral remainder therefore gives
\begin{equation}\label{eq:integralremainder}
 \begin{split}
 Y(t_{j+1})&=P_{j+1}(X_j)+R_j,\\
 R_j&=\int_0^{h_j}\frac{(h_j-s)^m}{m!}
                 F_m(Y(t_j+s))\dd s.
 \end{split}
\end{equation}
Using~\eqref{eq:remainderbound}, we obtain
\[
 \norm{R_j}
 \le\frac{(m+1)!\tau_j}{h_j^{m+1}}
       \int_0^{h_j}\frac{(h_j-s)^m}{m!}\dd s
 =\tau_j.
\]

Both $X_j$ and $z_j$ belong to the convex set $S_j$, so their
joining segment lies in $S_j$. The fundamental theorem of
calculus gives
\[
 P_{j+1}(X_j)-P_{j+1}(z_j)
 =\int_0^1DP_{j+1}\bigl(z_j+u(X_j-z_j)\bigr)
              (X_j-z_j)\dd u.
\]
It follows from~\eqref{eq:mapbound} that
\[
 \norm{P_{j+1}(X_j)-P_{j+1}(z_j)}
 \le A_j\norm{X_j-z_j}\le A_jE_j.
\]
Combining this inequality with~\eqref{eq:integralremainder},
\eqref{eq:defect}, and~\eqref{eq:recurrence}, we conclude that
\begin{align*}
 \norm{Y(t_{j+1})-z_{j+1}}
 &\le\norm{R_j}
       +\norm{P_{j+1}(X_j)-P_{j+1}(z_j)}
       +\norm{P_{j+1}(z_j)-z_{j+1}}\\
 &\le\tau_j+A_jE_j+\eta_j\le E_{j+1}.
\end{align*}
This proves the induction step. Hence the solution exists on
$[0,t_N]$ and satisfies all the stated error bounds.
The containment bounds follow from~\eqref{eq:displacementproof}
and~\eqref{eq:escape} at each step.
\end{proof}

\begin{remark}
\begin{itemize}
\item
In many numerical methods, it is natural to prescribe the step
sizes in advance. For the Round Taylor method, however, each
step must satisfy
\[
E_j+B_jh_j<\varepsilon_j,
\]
where the speed bound $B_j$ depends on the neighborhood determined
by $\varepsilon_j$. Thus, increasing $\varepsilon_j$ to accommodate
a prescribed step may also increase $B_j$.

A simple implementation fixes $\varepsilon_j=\varepsilon>0$.
At each step, one first computes $B_j$ on the corresponding
neighborhood and then selects $h_j$ to satisfy the strict
inequality. If $E_j\ge\varepsilon$, the procedure stops because
this inequality cannot hold.

\item
A natural practical choice is to use a fixed, sufficiently small
rounding grid $H_j=H$ for $j=0,\ldots,N$ and take $\eta_j=H$
at every step.
If different grids $H_j$ have been proposed, choosing
$H=\min_{0\le j\le N}H_j$ makes every prescribed rounding-error bound
no larger. This can improve the final error estimate with
a modest increase in computational cost.
We retain the separate quantities $H_j$ and $\eta_j$ in the
theorem for generality.

\end{itemize}
\end{remark}

\section{Polynomial vector fields}\label{sec:polynomial}

When the coordinates of $f$ are polynomials, the Taylor functions
$F_0,\ldots,F_m$, the predictors, and their derivatives are also
polynomials. In this section, we explain how this property allows
a computer to carry out the estimates required by
Theorem~\ref{thm:main} and produce an approximation with a
rigorous error bound that accounts for rounding at each step.

The application of an error estimate usually requires preliminary
bounds on the vector field and its derivatives. For polynomial
vector fields, these bounds can be computed as part of the
procedure, without deriving them by hand. With rational
coefficients, all steps and bounds after the initial rounding
can be computed using exact rational arithmetic.

We first explain how to bound a polynomial on a box. We then
give an algorithm that uses these bounds to select the step
sizes and compute the approximations and their error bounds.

\subsection{Polynomial bounds on a box}
\label{subsec:polynomialbounds}

\begin{lemma}[The range of a monomial]\label{lem:monomial}
Let $Q=\prod_{i=1}^n[a_i,b_i]$, with $a_i\le b_i$, and let
\[
 q(y)=c\prod_{i=1}^n y_i^{\alpha_i},\qquad c\in\R,\qquad
 \alpha_i\in\{0,1,2,\ldots\}.
\]
The exact range $[\ell_i,u_i]$ of the factor $y_i^{\alpha_i}$ is
given by the following rules:
\[
 (\ell_i,u_i)=
 \begin{cases}
 (1,1),&\alpha_i=0,\\
 (a_i^{\alpha_i},b_i^{\alpha_i}),&\alpha_i\text{ is odd},\\
 (0,\max\{a_i^{\alpha_i},b_i^{\alpha_i}\}),
     &\alpha_i>0\text{ is even and }a_i\le0\le b_i,\\
 (\min\{a_i^{\alpha_i},b_i^{\alpha_i}\},
  \max\{a_i^{\alpha_i},b_i^{\alpha_i}\}),
     &\alpha_i>0\text{ is even and }0\notin[a_i,b_i].
 \end{cases}
\]
A factor with exponent zero means the constant $1$, including
when its coordinate is zero. Set $L_0=U_0=1$. For $i=1,\ldots,n$,
let $L_i$ and $U_i$ be the minimum and maximum of the four numbers
\[
 L_{i-1}\ell_i,\quad L_{i-1}u_i,\quad
 U_{i-1}\ell_i,\quad U_{i-1}u_i.
\]
Then the exact minimum and maximum of $q$ on $Q$ are
\[
 (\min_Q q,\max_Q q)=
 \begin{cases}
 (cL_n,cU_n),&c\ge0,\\
 (cU_n,cL_n),&c<0.
 \end{cases}
\]
\end{lemma}
\begin{proof}
The power rules follow from monotonicity on each side of zero.
The minimum and maximum of a product $uv$ on a rectangle occur
among its four corner values. Since successive factors involve
different coordinates, their endpoint values can be chosen
independently, so the recursion gives exact extrema. Multiplication
by a negative coefficient reverses their order.
\end{proof}

\begin{lemma}[Bounds for a polynomial]\label{lem:polynomial}
Let $Q=\prod_{i=1}^n[a_i,b_i]$, with $a_i\le b_i$, and let
\[
 p(y)=\sum_\alpha c_\alpha y^\alpha,\qquad
 y^\alpha=\prod_{i=1}^n y_i^{\alpha_i},
\]
be a real polynomial after collecting like terms, where the sum
is finite and $\alpha_i\in\{0,1,2,\ldots\}$. For each term, let
\[
 \ell_\alpha=\min_{y\in Q}c_\alpha y^\alpha,\qquad
 u_\alpha=\max_{y\in Q}c_\alpha y^\alpha.
\]
Define
\begin{equation}\label{eq:polynomialbound}
 L_Q(p)=\sum_\alpha\ell_\alpha,\qquad
 U_Q(p)=\sum_\alpha u_\alpha,\qquad
 M_Q(p)=\max\{|L_Q(p)|,|U_Q(p)|\}.
\end{equation}
For the zero polynomial, the sums are empty and all three
quantities are zero. Then
\[
 L_Q(p)\le p(y)\le U_Q(p),\qquad |p(y)|\le M_Q(p)
 \quad(y\in Q).
\]
\end{lemma}
\begin{proof}
For every $y\in Q$, each term satisfies
$\ell_\alpha\le c_\alpha y^\alpha\le u_\alpha$.
Adding these inequalities gives the lower and upper bounds for
$p(y)$. The absolute-value bound follows from the definition of
$M_Q(p)$.
\end{proof}

\subsection{The algorithm}\label{subsec:algorithm}
Theorem~\ref{thm:main} allows prescribed time steps. To check
them, we still need regions on which the bounds hold. In
particular, $B_j$ depends on $\varepsilon_j$, so the condition
$E_j+B_jh_j<\varepsilon_j$ does not determine the radius before
the speed bound is known. The following procedure chooses the
radius first and then computes a sufficient step.

\paragraph{Inputs.}
The inputs are:
\begin{itemize}
\item a polynomial vector field $f=(f_1,\ldots,f_n):\R^n\to\R^n$
      with rational coefficients;
\item an initial state $y_0\in\R^n$ and a positive rational
      initial rounding grid $H_0$;
\item a positive integer Taylor order $m$ and a rational target
      time $T\ge0$;
\item a rational neighborhood radius $\varepsilon>0$ and a
      rational rounding grid $H>0$.
\end{itemize}
We use $\varepsilon_j=\varepsilon$ for $j=0,\ldots,N-1$
and $H_j=H$ for $j=1,\ldots,N$.
The steps $h_j$,
their number $N$, and all error and auxiliary bounds are computed
by the procedure. 

The initial state $y_0$ may have irrational coordinates; the
procedure requires access to its exact initial rounding
${\cal R}_{H_0}(y_0)$.

\paragraph{Possible outputs.}
Let $N\ge0$ be the number of completed steps. The procedure
returns one of the following:
\begin{enumerate}
\item \emph{Success:} $t_N=T$, a rational vector $z_N$, and
      a rational bound $E_N$, with the guarantees that the exact
      solution exists on $[0,T]$ and
      $\norm{Y(T)-z_N}\le E_N$.
\item \emph{Not certified at $T$:} a stopping reason and the
      last completed time $t_N<T$, together with $z_N,E_N$.
      The completed prefix remains certified: the solution
      exists on $[0,t_N]$ and
      $\norm{Y(t_i)-z_i}\le E_i$ for $i=0,\ldots,N$.
      This output does not imply that the solution fails to
      exist at $T$.
\end{enumerate}
In either case, the stored data include the points and errors
with indices $i=0,\ldots,N$, and the steps and their bounds
with indices $j=0,\ldots,N-1$. When $N=0$, only the initial
point and its error bound are stored.

\paragraph{Procedure.}
Set $z_0={\cal R}_{H_0}(y_0)$, $E_0=H_0$, and $t_0=0$.
For $T=0$, return success with these initial values. If $T>0$ and
$E_0\ge\varepsilon$, return not certified. Otherwise, compute
the polynomials $F_0,\ldots,F_m$ once. Starting at $j=0$,
perform the following steps.
\begin{enumerate}
\item If $t_j=T$, return success. Otherwise, return not
      certified if $E_j\ge\varepsilon$ or an imposed step
      limit has been reached.
\item Form the two boxes
      \[
       \begin{split}
       \overline{C_j}&=\prod_{i=1}^n
           [(z_j)_i-\varepsilon,(z_j)_i+\varepsilon],\\
       S_j&=\prod_{i=1}^n[(z_j)_i-E_j,(z_j)_i+E_j].
       \end{split}
      \]
      Using Lemma~\ref{lem:polynomial}, compute
      \begin{equation}\label{eq:polynomialspeed}
       B_j=\max_{1\le i\le n}M_{\overline{C_j}}(f_i).
      \end{equation}
      If $B_j>0$, choose
      \begin{equation}\label{eq:adaptivestep}
       h_j=\min\left\{T-t_j,\;
          H\left\lfloor\frac{\varepsilon-E_j}{2B_jH}\right\rfloor
                   \right\}.
      \end{equation}
      If $B_j=0$, take $h_j=T-t_j$. If $h_j=0$, return not
      certified, with the reason that the step rule produced
      no positive step.
\item With $h_j$ fixed, form the predictor $P_{j+1}$
      from~\eqref{eq:predictor}. Write $(F_m)_i$ and
      $(P_{j+1})_i$ for the $i$th coordinate polynomials and
      compute
      \begin{equation}\label{eq:polynomialconstants}
       \begin{split}
       \tau_j&=\frac{h_j^{m+1}}{(m+1)!}
          \max_{1\le i\le n}M_{\overline{C_j}}((F_m)_i),\\
       A_j&=\max_{1\le i\le n}\sum_{k=1}^n
          M_{S_j}\!\left(
             \frac{\partial(P_{j+1})_i}{\partial y_k}\right).
       \end{split}
      \end{equation}
      The last formula bounds the maximum absolute row sum of
      the complete predictor derivative. Its polynomial entries
      are collected before applying the monomial bounds.
\item Compute and store
      \begin{equation}\label{eq:adaptiveupdate}
       \begin{split}
       z_{j+1}&={\cal R}_H(P_{j+1}(z_j)),\\
       \widehat E_{j+1}&=A_jE_j+\tau_j+H,\\
       E_{j+1}&=H\left(1+\left\lfloor
                         \widehat E_{j+1}/H\right\rfloor\right),\\
       t_{j+1}&=t_j+h_j.
       \end{split}
      \end{equation}
      Continue with index $j+1$.
\end{enumerate}
All these calculations use rational arithmetic. If a computation
limit interrupts an unfinished step, return not certified and
retain only the fully completed steps.

For $B_j>0$, the step rule gives
\[
 E_j+B_jh_j\le\frac{\varepsilon+E_j}{2}<\varepsilon.
\]
For $B_j=0$, the same strict containment follows from
$E_j<\varepsilon$. The point defect is less than $H$, so
we take $\eta_j=H$, and the upward error rounding satisfies
\[
 \widehat E_{j+1}<E_{j+1}\le\widehat E_{j+1}+H.
\]
Thus every completed step satisfies Theorem~\ref{thm:main}.
For $T=0$, the initial rounding bound gives the output bound directly.
The final error $E_N$ may be at least $\varepsilon$ without
invalidating the result, because success is checked before the
condition for attempting another step.

The last step is shortened to $T-t_j$ before its predictor and
bounds are computed. Every positive step that does not reach
$T$ is at least $H$, so the procedure terminates after finitely
many steps.

An implementation of this fixed-radius, fixed-grid procedure is
available at
\begin{center}
 \url{https://www.oscarperdomo.org/round-taylor}.
\end{center}
It uses exact rational arithmetic and displays the stored points,
errors, and step bounds. It also exports certificates and provides
a separate Python checker. The application uses automatic initial
rounding. The current application accepts up to
six variables and orders $1\le m\le10$, subject to its stated
computation limits.

\section{Examples that can be checked by hand}\label{sec:examples}
We first illustrate the polynomial bounds and then carry out
the procedure by hand in one and two variables. In these examples,
we take $H_0=H$, choose $z_0={\cal R}_H(y_0)$, and set $E_0=H_0$.

\subsection{Bounds on an interval}
Take $p(y)=2y-4y^2$ and $Q=[3,5]$. The ranges of the two
terms are $[6,10]$ and $[-100,-36]$. Adding their endpoints
gives
\[
 L_Q(p)=-94,\qquad U_Q(p)=-26,\qquad M_Q(p)=94.
\]
Consequently $-94\le p(y)\le-26$ and $|p(y)|\le94$ on $Q$.
These are bounds for the polynomial; its different terms need
not attain their extrema at the same point.

\subsection{A scalar equation}
\paragraph{Inputs.}
Take
\[
 \begin{gathered}
 f(y)=y^2,\qquad y_0=1,\qquad m=1,\\
 T=\frac1{10},\qquad \varepsilon=\frac12,\qquad H=\frac1{100}.
 \end{gathered}
\]
\paragraph{Initialization.}
The initial stored data and Taylor functions are
\[
 \begin{gathered}
 z_0=1,\qquad E_0=\frac1{100},\qquad t_0=0,\\
 F_0(y)=y^2,\qquad F_1(y)=2y^3.
 \end{gathered}
\]
Start with $j=0$.
\begin{enumerate}
\item Since $t_0<T$ and $E_0<\varepsilon$, continue.
\item Form the intervals
      \[
       \overline{C_0}=[1/2,3/2],\qquad
       S_0=[99/100,101/100].
      \]
      The polynomial bound gives $B_0=9/4$. Hence
      \[
       h_0=\min\left\{\frac1{10},\frac1{100}
                      \left\lfloor\frac{98}{9}\right\rfloor
                \right\}=\frac1{10}.
      \]
      The containment inequality is
      $E_0+B_0h_0=47/200<1/2$.
\item With this step, the predictor is
      $P_1(y)=y+y^2/10$. On $S_0$, its derivative is
      $1+y/5$, so $A_0=601/500$. On $\overline{C_0}$,
      $M_{\overline{C_0}}(F_1)=27/4$, giving
      \[
       \tau_0=\frac{h_0^2}{2}\frac{27}{4}=\frac{27}{800}.
      \]
\item Round the predictor value and update the error and time:
      \[
       \begin{split}
       z_1&={\cal R}_{1/100}(P_1(1))=\frac{11}{10},\\
       \widehat E_1&=\frac{601}{500}\frac1{100}
                    +\frac{27}{800}+\frac1{100}
                    =\frac{5577}{100000},\\
       E_1&=\frac1{100}\left(1+
                    \left\lfloor\frac{5577}{1000}\right\rfloor
                        \right)=\frac3{50},\qquad t_1=\frac1{10}.
       \end{split}
      \]
\end{enumerate}
\paragraph{Output.}
At the next stopping check, $t_1=T$, so the procedure returns
\emph{Success}, with $N=1$, $z_1=11/10$, and $E_1=3/50$.
Thus the exact solution exists on $[0,1/10]$ and
$|Y(1/10)-11/10|\le3/50$. The exact solution
$Y(t)=1/(1-t)$, for $0\le t<1$, gives the independent check
$|Y(1/10)-11/10|=1/90<3/50$.

If the same $f,y_0,m,T$ are used with
$\varepsilon=H=1/100$, the procedure instead returns
\emph{Not certified at $T$} immediately after initialization,
because $E_0=H\ge\varepsilon$. Its saved output is
$N=0$, $t_0=0$, $z_0=1$, and $E_0=1/100$.
The stopping reason concerns these parameter choices; the exact
solution still exists at $T$.

\subsection{A coupled system in two variables}\label{subsec:coupledexample}
\paragraph{Inputs.}
Take
\[
 \begin{gathered}
 f(u,v)=(v,uv),\qquad y_0=(1,1),\qquad m=1,\\
 T=\frac9{50},\qquad \varepsilon=\frac12,\qquad H=\frac1{100}.
 \end{gathered}
\]
All terminating decimals below denote exact rational numbers.
\paragraph{Initialization.}
We obtain
\[
 \begin{gathered}
 z_0=(1,1),\qquad E_0=0.01,\qquad t_0=0,\\
 F_1(u,v)=(uv,v^2+u^2v).
 \end{gathered}
\]
\paragraph{First iteration, $j=0$.}
\begin{enumerate}
\item Since $t_0<T$ and $E_0<\varepsilon$, continue.
\item Form
      \[
       \overline{C_0}=[0.5,1.5]^2,\qquad S_0=[0.99,1.01]^2.
      \]
      The coordinate bounds for $f$ are $1.5$ and $2.25$,
      so $B_0=2.25$. The step rule gives
      \[
       h_0=\min\{0.18,\;0.01\lfloor98/9\rfloor\}=0.10,
       \qquad E_0+B_0h_0=0.235<0.5.
      \]
\item Form the predictor and its derivative:
      \[
       P_1(u,v)=(u+0.10v,\ v+0.10uv),\qquad
       DP_1(u,v)=
          \begin{pmatrix}1&0.10\\0.10v&1+0.10u\end{pmatrix}.
      \]
      On $S_0$, the row-sum bounds are $1.10$ and
      $1+0.10(1.01+1.01)=1.202$, so $A_0=1.202$.
      On $\overline{C_0}$, the coordinate bounds for $F_1$
      are $1.5^2=2.25$ and $1.5^2+1.5^3=5.625$. Therefore
      \[
       \tau_0=\frac{0.10^2}{2}(5.625)=0.028125.
      \]
\item Since $P_1(z_0)=(1.10,1.10)$ is already on the grid,
      \[
       \begin{split}
       z_1&=(1.10,1.10),\\
       \widehat E_1&=1.202(0.01)+0.028125+0.01=0.050145,\\
       E_1&=0.01(1+\lfloor5.0145\rfloor)=0.06,
       \qquad t_1=0.10.
       \end{split}
      \]
\end{enumerate}
\paragraph{Second iteration, $j=1$.}
\begin{enumerate}
\item Since $t_1<T$ and $E_1<\varepsilon$, continue.
\item The new boxes are
      \[
       \overline{C_1}=[0.6,1.6]^2,\qquad S_1=[1.04,1.16]^2.
      \]
      The speed bound is $B_1=1.6^2=2.56$. Since $T-t_1=0.08$,
      \[
       h_1=\min\{0.08,\;0.01\lfloor275/32\rfloor\}=0.08,
       \qquad E_1+B_1h_1=0.2648<0.5.
      \]
\item Now $P_2(u,v)=(u+0.08v,\ v+0.08uv)$.
      The row-sum bounds on $S_1$ are $1.08$ and
      $1+0.08(1.16+1.16)=1.1856$, so $A_1=1.1856$.
      The largest coordinate bound for $F_1$ on
      $\overline{C_1}$ is $1.6^2+1.6^3=6.656$. Hence
      \[
       \tau_1=\frac{0.08^2}{2}(6.656)=0.0212992.
      \]
\item Evaluate $P_2(z_1)=(1.188,1.1968)$ and round each
      coordinate downward to the grid. This gives
      \[
       \begin{split}
       z_2&=(1.18,1.19),\\
       \widehat E_2&=1.1856(0.06)+0.0212992+0.01=0.1024352,\\
       E_2&=0.01(1+\lfloor10.24352\rfloor)=0.11,
       \qquad t_2=0.18.
       \end{split}
      \]
\end{enumerate}
\paragraph{Output.}
At the next stopping check, $t_2=T$. The procedure returns
\emph{Success}, with
\[
 N=2,\qquad z_2=(59/50,119/100),\qquad E_2=11/100.
\]
Thus $Y=(u,v)$ exists on $[0,9/50]$ and
\[
 \left\|Y(9/50)-(59/50,119/100)\right\|_\infty\le\frac{11}{100}.
\]
To enter this example in the application, use the field
\texttt{\{x2, x1*x2\}} and initial vector \texttt{\{1,1\}},
with the rational inputs $T=9/50$, $m=1$,
$\varepsilon=1/2$, and $H=1/100$.

\section{Computer examples}\label{sec:computerexamples}
The calculations in this section are done by computer, using the
polynomial procedure in Section~\ref{subsec:algorithm}. We state
the initial data, $m$, $\varepsilon$, $H$, the target time, and
the resulting error. Both examples use the fixed grid $H=10^{-50}$
and take $H_0=H$, $z_0={\cal R}_H(y_0)$, $E_0=H_0$. The steps,
rounded points, and error bounds are computed using exact rational
arithmetic. The accompanying scripts carry out these computations.

\subsection{The scalar equation}
Consider the initial-value problem
\[
 Y'(t)=Y(t)^2,\qquad Y(0)=1.
\]
Use
\[
 \boxed{m=10,\qquad \varepsilon=\frac1{50},\qquad H=10^{-50},
        \qquad T=\frac12.}
\]
The polynomial procedure reaches $T$ in $103$ steps and gives
\[
 \boxed{|Y(1/2)-z_{103}|\le E_{103}<3.3\times10^{-21}.}
\]
This bounds the absolute error in the stored approximation at
$T=1/2$, including rounding at each step.

Here $F_\ell(y)=(\ell+1)!y^{\ell+2}$, so the polynomial bounds
used by the computer take the simple forms
\[
 \begin{aligned}
 B_j&=(z_j+\varepsilon)^2,\\
 \tau_j&=h_j^{m+1}(z_j+\varepsilon)^{m+2},\\
 A_j&=\sum_{k=0}^m(k+1)h_j^k(z_j+E_j)^k.
 \end{aligned}
\]
Every interval in this run has positive endpoints. The steps,
rounded points, and error bounds are computed by the same formulas
as in the hand examples. The explicit solution $Y(t)=1/(1-t)$
provides an exact check of the error at every stored time.

\subsection{A four-body problem}
\label{sec:fourbody}

The polynomial procedure can also be useful when the original
differential equation contains radicals. In this section we explain
how to introduce additional variables that remove the radicals from
a gravitational four-body problem. The resulting system has the
polynomial form required in Section~\ref{subsec:algorithm}.
Our aim is to assess how close the proposed motion is to being periodic.

\subsubsection{The differential equations}
\begin{figure}[htbp]
\centering
\begin{tikzpicture}[>=Stealth, font=\small,
 bodyone/.style={circle,draw=inkblue,fill=inkblue,inner sep=2.4pt},
 bodymu/.style={circle,draw=leaf,fill=leaf,inner sep=2.4pt}]
\begin{scope}[xshift=-3.2cm]
 \node[font=\small\bfseries] at (0,2.35) {Opposite pairs at any time};
 \draw[->,gray!65] (-2.65,0)--(2.65,0);
 \draw[->,gray!65] (0,-1.95)--(0,1.95);
 \draw[densely dashed,inkblue!60] (-1.9,-.95)--(1.9,.95);
 \draw[densely dashed,leaf!70] (.8,-1.4)--(-.8,1.4);
 \node[bodyone,label={[align=center]right:{$\mathbf u$\\mass $1$}}] at (1.9,.95) {};
 \node[bodyone,label={[align=center]left:{$-\mathbf u$\\mass $1$}}] at (-1.9,-.95) {};
 \node[bodymu,label={[align=center]above:{$\mathbf v$; mass $\mu$}}] at (-.8,1.4) {};
 \node[bodymu,label={[align=center]below:{$-\mathbf v$; mass $\mu$}}] at (.8,-1.4) {};
 \draw[->,thick] (.53,.28) arc[start angle=28,end angle=208,radius=.60];
 \node at (-.85,.38) {$\pi$};
 \node[below right,inner sep=2pt] at (0,0) {$O$};
 \node[align=center,font=\footnotesize] at (0,-2.65)
  {A half-turn exchanges the two bodies\\within each equal-mass pair.};
\end{scope}
\begin{scope}[xshift=3.2cm]
 \node[font=\small\bfseries] at (0,2.35) {Initial positions and velocities};
 \draw[->,gray!65] (-2.65,0)--(2.65,0) node[right] {$x$};
 \draw[->,gray!65] (0,-1.95)--(0,1.95) node[above] {$y$};
 \node[bodyone] at (1.95,0) {};
 \node[bodyone] at (-1.95,0) {};
 \node[bodymu] at (0,1.25) {};
 \node[bodymu] at (0,-1.25) {};
 \node[below right] at (1.95,0) {$\mathbf u(0)$};
 \node[above left] at (-1.95,0) {$-\mathbf u(0)$};
 \node[above left] at (0,1.25) {$\mathbf v(0)$};
 \node[below right] at (0,-1.25) {$-\mathbf v(0)$};
 \draw[->,very thick,inkblue] (1.95,.10)--(1.95,.95) node[right] {$\mathbf p(0)$};
 \draw[->,very thick,inkblue] (-1.95,-.10)--(-1.95,-.95) node[left] {$-\mathbf p(0)$};
 \draw[->,very thick,leaf] (.10,1.25)--(1.05,1.25) node[above] {$\mathbf q(0)$};
 \draw[->,very thick,leaf] (-.10,-1.25)--(-1.05,-1.25) node[below] {$-\mathbf q(0)$};
 \node[below right,inner sep=2pt] at (0,0) {$O$};
 \node[align=center,font=\footnotesize] at (0,-2.65)
  {Reflection in the $x$-axis, together with time\\reversal, preserves these initial data\\after exchanging the two mass-$\mu$ bodies.};
\end{scope}
\end{tikzpicture}
\caption{The symmetries of the four-body configuration. The half-turn
symmetry on the left allows us to describe all four positions using
only $\mathbf u$ and $\mathbf v$. The arrows on the right are the
initial velocities for the signs of the angular velocities in
\eqref{eq:fourbodydata}. The diagrams are schematic and are not orbit plots.}
\label{fig:fourbodysymmetry}
\end{figure}
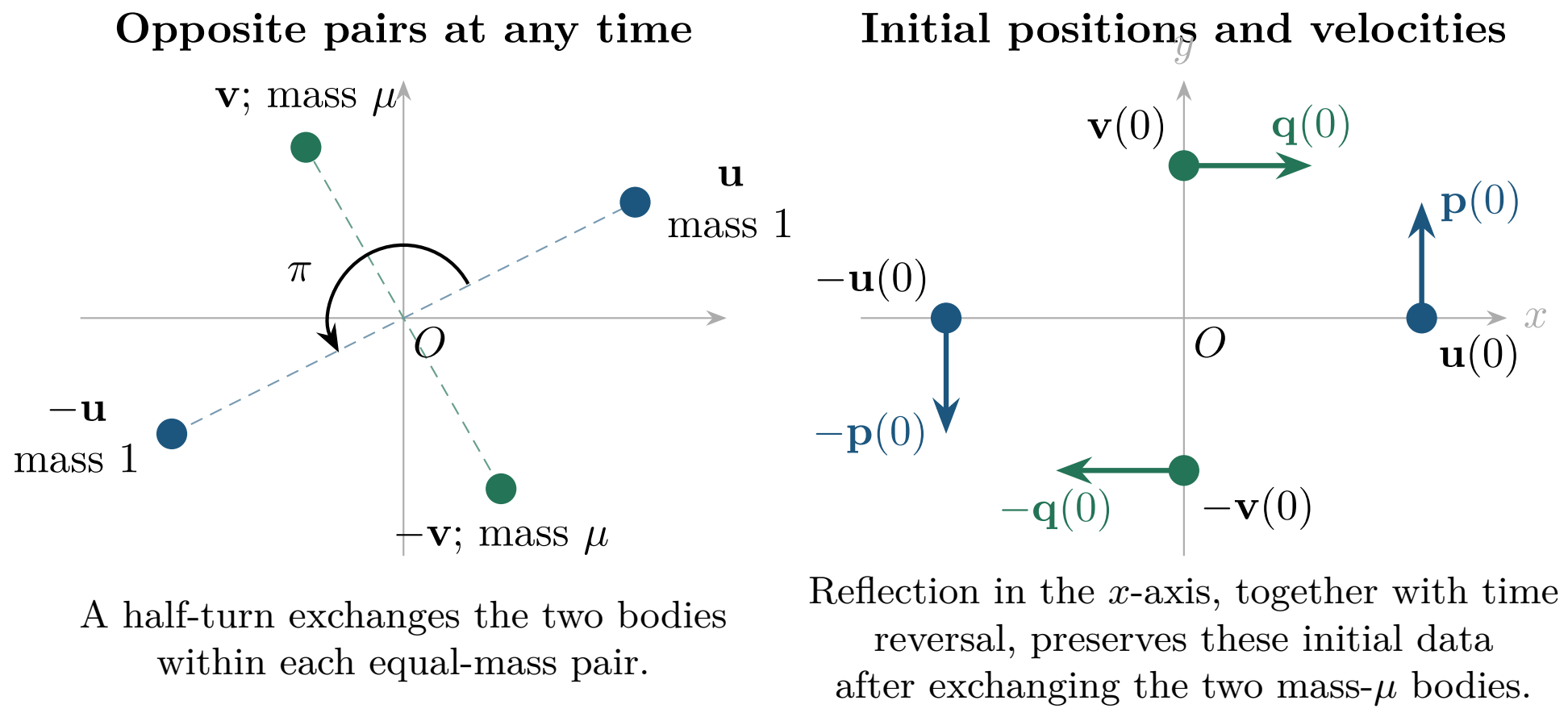

We consider the four-body configuration in
\cite[Section~2, equations~(13)--(14)]{Perdomo2026NBody}.
The gravitational constant is $1$, the first opposite pair has masses
$1,1$, and the second has masses $\mu,\mu$, where $\mu>0$ denotes
the mass called $m_2$ in that paper. Their positions are
\[
 \mathbf u,-\mathbf u,\mathbf v,-\mathbf v,\qquad
 \mathbf u=r_1(\cos\theta,\sin\theta),\quad
 \mathbf v=r_2(\cos(\beta+\pi/2),\sin(\beta+\pi/2)).
\]
The two bodies in each opposite pair have the same mass. Rotation
through $\pi$ exchanges them and preserves Newton's equations.
By uniqueness, initial positions and velocities in opposite pairs
remain in opposite pairs throughout the collision-free solution.
This is the symmetry used to reduce the four-body system to the
two vectors $\mathbf u$ and $\mathbf v$; see
Figure~\ref{fig:fourbodysymmetry}.

Writing $\delta=\theta-\beta$ and
\[
 d_1=\sqrt{r_1^2+r_2^2-2r_1r_2\sin\delta},\qquad
 d_2=\sqrt{r_1^2+r_2^2+2r_1r_2\sin\delta},
\]
the differential equations from that manuscript are
\begin{equation}\label{eq:fourbodypolar}
\begin{aligned}
 \ddot r_1&=r_1\dot\theta^2-\frac{1}{4r_1^2}
   +\mu\left(\frac{r_2\sin\delta-r_1}{d_1^3}
             +\frac{-r_2\sin\delta-r_1}{d_2^3}\right),\\
 \ddot\theta&=\frac{\mu r_2\cos\delta}{r_1}
                \left(\frac1{d_1^3}-\frac1{d_2^3}\right)
                -\frac{2\dot r_1\dot\theta}{r_1},\\
 \ddot r_2&=r_2\dot\beta^2-\frac{\mu}{4r_2^2}
   +\frac{r_1\sin\delta-r_2}{d_1^3}
   +\frac{-r_1\sin\delta-r_2}{d_2^3},\\
 \ddot\beta&=\frac{r_1\cos\delta}{r_2}
                \left(\frac1{d_2^3}-\frac1{d_1^3}\right)
                -\frac{2\dot r_2\dot\beta}{r_2}.
\end{aligned}
\end{equation}
This vector field is not polynomial. The distances contain square
roots, and the polar equations also contain trigonometric functions
and divisions by the radii. Cartesian positions remove the
trigonometric functions. Reciprocal-distance variables remove the
remaining divisions and radicals from the vector field.

\subsubsection{Removing the radicals}
In this section $|\mathbf w|=\sqrt{w_1^2+w_2^2}$ denotes the
Euclidean length of a planar vector, and the dot denotes the
Euclidean inner product. Set
\[
 \mathbf p=\dot{\mathbf u},\qquad
 \mathbf q=\dot{\mathbf v},\qquad
 a=\frac1{|\mathbf u|},\quad b=\frac1{|\mathbf v|},\quad
 c=\frac1{|\mathbf u-\mathbf v|},\quad
 d=\frac1{|\mathbf u+\mathbf v|}.
\]
For example, the force term
$(\mathbf v-\mathbf u)/|\mathbf v-\mathbf u|^3$ becomes
$c^3(\mathbf v-\mathbf u)$. Differentiating the definitions of
$a,b,c,d$ gives differential equations for these new variables.
The resulting first-order system is
\begin{equation}\label{eq:fourbodypolynomial}
\begin{aligned}
 \dot{\mathbf u}&=\mathbf p,&
 \dot{\mathbf v}&=\mathbf q,\\
 \dot{\mathbf p}&=-\frac{a^3}{4}\mathbf u
       +\mu\bigl(c^3(\mathbf v-\mathbf u)
                    -d^3(\mathbf u+\mathbf v)\bigr),\\
 \dot{\mathbf q}&=-\frac{\mu b^3}{4}\mathbf v
       +c^3(\mathbf u-\mathbf v)-d^3(\mathbf u+\mathbf v),\\
 \dot a&=-a^3\mathbf u\cdot\mathbf p,&
 \dot b&=-b^3\mathbf v\cdot\mathbf q,\\
 \dot c&=-c^3(\mathbf u-\mathbf v)\cdot(\mathbf p-\mathbf q),\\
 \dot d&=-d^3(\mathbf u+\mathbf v)\cdot(\mathbf p+\mathbf q).
\end{aligned}
\end{equation}
The state is $Y=(\mathbf u,\mathbf v,\mathbf p,\mathbf q,a,b,c,d)$.
It has twelve scalar coordinates: two for each of the four vectors
and one for each of $a,b,c,d$. Every right-hand side is polynomial. Its coefficients
are rational when $\mu$ is rational. The factors $1/4$ in the
self-pair forces follow from the distance $2|\mathbf u|$ or
$2|\mathbf v|$ between opposite bodies.

The two systems describe the same motion with the initial data
below. For a reciprocal variable $s$ corresponding to a vector
$\mathbf w$, set $C=|\mathbf w|^2s^2-1$.
The last four equations give
\[
 C'=-2s^2(\mathbf w\cdot\mathbf w')C.
\]
Since $C(0)=0$, we have $C=0$ throughout the solution.
The reciprocal identities are therefore preserved, and the
physical distances remain positive while the polynomial solution
exists.

\subsubsection{Initial data}
\begin{figure}[!ht]
\centering
\input{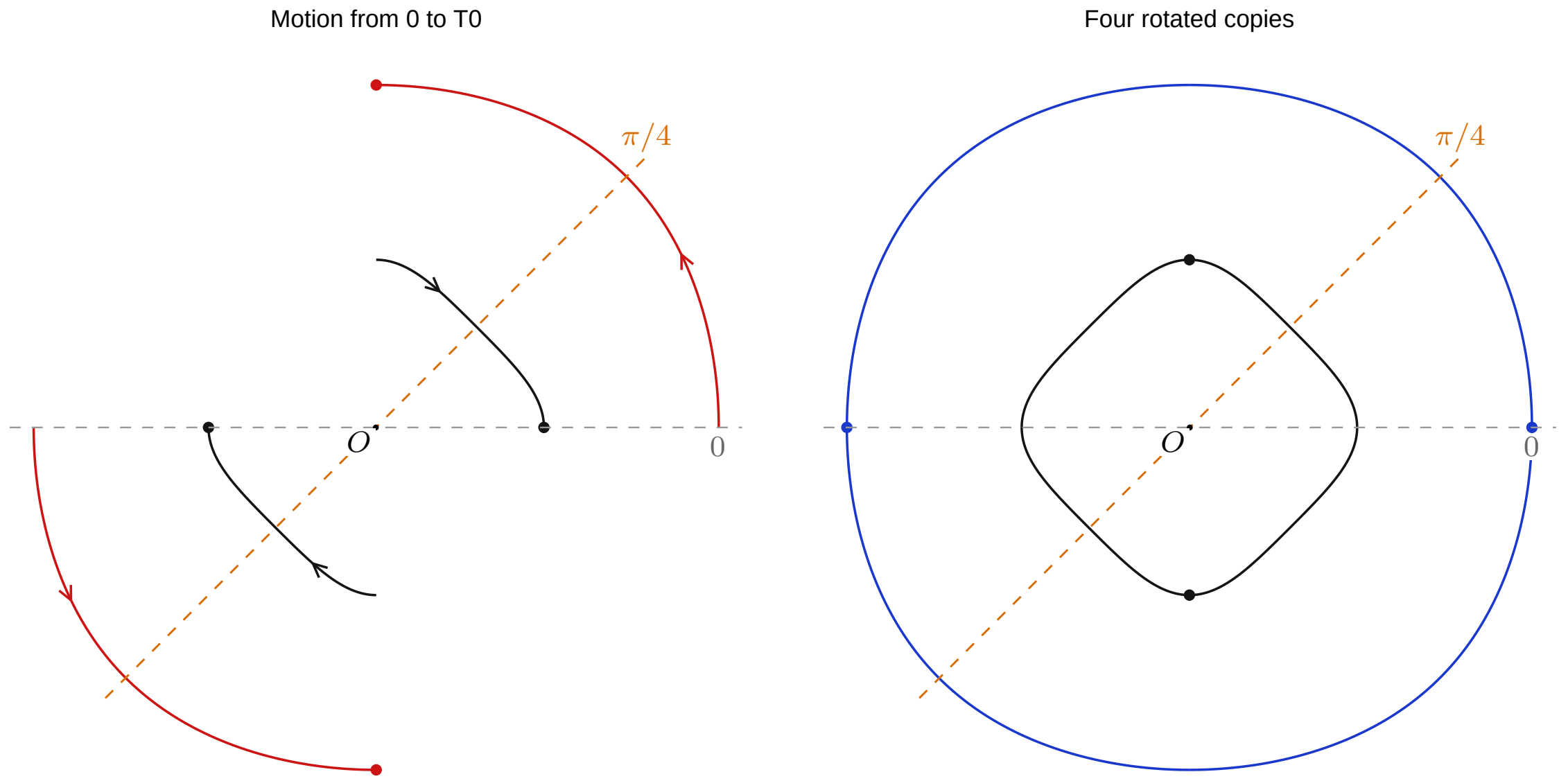}
\caption{Numerical drawings of the four-body motion for the
$90^\circ$ data in~\eqref{eq:fourbodydata}. Left: the four
paths from $0$ to $T_0$; red corresponds to the mass-$1$ pair and
black to the mass-$\mu$ pair. The dots mark the positions at
$T_0$. Right: the pattern formed by four rotated copies of the
first paths, with blue for the mass-$1$ pair and black for the
mass-$\mu$ pair. The dots in the right panel mark the initial
positions. These drawings illustrate the proposed periodic
motion. The dashed lines through $O$ mark the initial reflection
axis (horizontal, gray) and the required midpoint reflection axis
at $\pi/4$ (orange). The computation below certifies the state and
the midpoint symmetry residuals at $t_*=T_0/2$.}
\label{fig:fourbodytrajectories}
\end{figure}

We use the $90^\circ$ row of Table~1 in~\cite{Perdomo2026NBody}.
These data are normalized so that $r_1(0)=1$:
\begin{equation}\label{eq:fourbodydata}
\begin{aligned}
 r_{1,0}&=1,&
 r_{2,0}&=0.4896144338645078,\\
 \omega_{1,0}&=1.8084629900120226,&
 \omega_{2,0}&=-1.2727625998400442,\\
 \mu&=1.5658646583060000,&
 T_0&=0.8472253023359529.
\end{aligned}
\end{equation}
Here $T_0$ is the proposed return time from that table. We use it
for the trajectory drawing in
Figure~\ref{fig:fourbodytrajectories}.

Every decimal in~\eqref{eq:fourbodydata} is interpreted as an
exact rational number. The initial physical data are
\[
\begin{aligned}
 \mathbf u(0)&=(r_{1,0},0),&
 \mathbf v(0)&=(0,r_{2,0}),\\
 \mathbf p(0)&=(0,r_{1,0}\omega_{1,0}),&
 \mathbf q(0)&=(-r_{2,0}\omega_{2,0},0).
\end{aligned}
\]
The remaining initial coordinates are
\[
 a(0)=\frac1{r_{1,0}},\qquad b(0)=\frac1{r_{2,0}},\qquad
 c(0)=d(0)=\frac1{\sqrt{r_{1,0}^2+r_{2,0}^2}}.
\]
For this computation, take $H_0=H$, choose
$z_0=\mathcal R_H(Y(0))$, and set $E_0=H_0$. The method includes this
initial error in its estimate, including the two irrational
reciprocal distances.

The initial velocities are perpendicular to their position vectors.
Reflection in the $x$-axis together with time reversal preserves
these data after exchanging the two mass-$\mu$ bodies, as shown
in Figure~\ref{fig:fourbodysymmetry}.

\subsubsection{The midpoint equations and why they suffice}
The goal is to check the four midpoint equations in the revised
manuscript at
\[
 t_*=\frac{T_0}{2}=0.42361265116797645.
\]
They are
\begin{equation}\label{eq:fourbodymidpoint}
 \dot r_1(t_*)=\dot r_2(t_*)=0,\qquad
 \delta(t_*)=\frac\pi2,\qquad \theta(t_*)=\frac\pi4.
\end{equation}
Here $H$ denotes the rounding grid and $t_*$ denotes the midpoint
time. Since $\delta=\theta-\beta$, the angular equations give
$\beta(t_*)=-\pi/4$. The polar angle of $\mathbf v$ is
$\beta+\pi/2$, so both position vectors lie on the line making
angle $\pi/4$ with the $x$-axis. The zero radial velocities mean
that both velocity vectors are perpendicular to this line.

Let $S$ be reflection in this line and let $R$ be rotation
through $\pi/2$. If~\eqref{eq:fourbodymidpoint} holds exactly,
reflection fixes both positions and negates both velocities
at $t_*$. Newton's equations are unchanged by reflection and
time reversal. Consequently the reflected motion, run backward
from $t_*$, has the same complete state there as the original
motion. Uniqueness therefore gives
\[
 \mathbf u(t_*+s)=S\mathbf u(t_*-s),\qquad
 \mathbf v(t_*+s)=S\mathbf v(t_*-s)
 \quad(0\le s\le t_*).
\]
The reflected half has the same interbody distances, so it is
also collision-free. At $s=t_*$, the initial positions and
velocities give
\begin{equation}\label{eq:fourbodyreturn}
 \begin{aligned}
 \mathbf u(T_0)&=R\mathbf u(0),&
 \mathbf p(T_0)&=R\mathbf p(0),\\
 \mathbf v(T_0)&=-R\mathbf v(0),&
 \mathbf q(T_0)&=-R\mathbf q(0).
 \end{aligned}
\end{equation}
The minus sign exchanges the two bodies in the second pair.
Rotation and this exchange preserve the equations and the
masses, so uniqueness allows the same motion to be continued
through successive copies. After four copies, the rotation and
the exchange are both the identity. Thus $4T_0$ is a period of
the labeled motion, though it need not be the smallest period.
This is the midpoint argument of
\cite[Lemma~2.2 and the return argument in Section~1]{Perdomo2026NBody}.

\subsubsection{The computer calculation}
We use the initial data above and choose
\[
 \boxed{m=4,\qquad \varepsilon=\frac1{625},\qquad H=10^{-50},
        \qquad T=t_*=\frac{T_0}{2}.}
\]
The polynomial procedure in Section~\ref{subsec:algorithm}
reaches $t_*=0.42361265116797645$ in $1774$ steps, with
\begin{equation}\label{eq:fourbodymidpointerror}
 \boxed{\|Y(t_*)-z_{1774}\|_\infty
       \le E_{1774}<5.9\times10^{-9}.}
\end{equation}
The last step is shortened to end exactly at $t_*$. This error
bounds each of the twelve coordinates and includes the rounding
at every step.

\begin{samepage}
The first eight coordinates of $z_{1774}$ give the following stored
position and velocity vectors at $t_*$. The bars denote stored
approximations. Each coordinate is written with exactly $50$
decimal places because $H=10^{-50}$; its certified error is still
bounded by~\eqref{eq:fourbodymidpointerror}.
\[
 \begin{aligned}
 \bar{\mathbf u}&=\begin{pmatrix}
  0.73141004583726970560542635335193881943717163805733\\
  0.73141004583730087236666405558947300876709723705726
 \end{pmatrix},\\[4pt]
 \bar{\mathbf v}&=\begin{pmatrix}
  0.29326824605887585305197515245893773326061930958294\\
  0.29326824605884877682147517874661312620365649026621
 \end{pmatrix},\\[4pt]
 \bar{\mathbf p}&=\begin{pmatrix}
  -1.40550721088382516253679578993486789208628104182917\\
  1.40550721088354199262816299075150820204509480403662
 \end{pmatrix},\\[4pt]
 \bar{\mathbf q}&=\begin{pmatrix}
  0.78971304083993280741800633128329718643883300826234\\
  -0.78971304083979975955668271879836334819145059753401
 \end{pmatrix}.
 \end{aligned}
\]
\end{samepage}

We recover the quantities in the midpoint equations directly
from the output:
\begin{equation}\label{eq:fourbodyrecover}
 \begin{aligned}
 \dot r_1&=a(u_xp_x+u_yp_y),&
 \dot r_2&=b(v_xq_x+v_yq_y),\\
 \theta&=\arctan\frac{u_y}{u_x},&
 \delta&=\theta+\arctan\frac{v_x}{v_y}.
 \end{aligned}
\end{equation}
The certified boxes keep $u_x>0$ and $v_y>0$ throughout the
calculation. Thus these inverse tangents use the continuous
angle branches with $\theta(0)=\beta(0)=0$ and
$\delta=\theta-\beta$.

Substituting the stored point in~\eqref{eq:fourbodyrecover}
gives the middle column below. Using its error enclosure gives
the last column, which bounds the absolute value of each
quantity for the exact solution. All quantities are evaluated
at $t_*$, and the angular errors are in radians.
\begin{center}
\begin{tabular}{@{}lll@{}}
\toprule
Quantity & Value from the stored point & Certified absolute bound\\
\midrule
$\dot r_1$ & $-1.5788\times10^{-13}$ & $<2.4\times10^{-8}$\\
$\dot r_2$ & $\phantom{-}1.4563\times10^{-13}$ & $<3.1\times10^{-8}$\\
$\delta-\pi/2$ & $\phantom{-}6.7469\times10^{-14}$ & $<2.8\times10^{-8}$\\
$\theta-\pi/4$ & $\phantom{-}2.1306\times10^{-14}$ & $<8\times10^{-9}$\\
\bottomrule
\end{tabular}
\end{center}
For the radial velocities, we insert the coordinate enclosure
in the first two formulas of~\eqref{eq:fourbodyrecover}.
For the angles, the exact differences can be written as
\[
 \begin{aligned}
 \theta-\frac\pi4
 &=\arctan\frac{u_y-u_x}{u_x+u_y},\\
 \delta-\frac\pi2
 &=\arctan\frac{u_yv_x-u_xv_y}{u_xv_x+u_yv_y}.
 \end{aligned}
\]
The position coordinates are positive on the endpoint enclosure,
so the denominators are positive and these identities use the
stated branches. The elementary inequality
$|\arctan s|\le |s|$ gives the two angular bounds from the same
coordinate enclosure. The accompanying code carries out these
evaluations using exact rational comparisons.

These estimates show how accurately the four midpoint equations
are satisfied by the solution with the printed initial data.
The symmetry argument applies when their residuals are exactly
zero. Proving that nearby parameters make them vanish would
require a separate root argument.
\FloatBarrier

\end{document}